# REPRESENTATION OF INTEGERS BY SUM OF THREE CUBES – A NEW APPROACH BASED ON SEED EQUATION

**Narinder Kumar Wadhawan**[1*] and **Priyanka Wadhawan**[2]

[1*]Civil Servant, Indian, Administrative Service Retired, B.Sc.(Engg.) From PEC, Chandigarh, India, Email: narinderkw@gmail.com

[2]Department Of Computer Sciences, Thapar Institute Of Engineering And Technology, Patiala, India, Email: priyanka.wadhawan@gmail.com

**Abstract**: We have proved in this paper that numbers can be expressed in algebraic form using one variable and two real rational quantities and thus sum of three cubes can also be expressed in algebraic form as a cubic polynomial. Using skeletal or seed equation, the polynomial can be transformed into a quadratic equation. A seed equation denotes a simple equation that represents a given integer as sum of three cubes including 0, for example, a seed equation for integer 2 is 1 cubed plus 1 cubed plus 0 cubed. Resultant quadratic equation can further be transformed into a linear equation which yields value of the variable and substitution of the value of the variable into the algebraic form of numbers results in the required solution. Notwithstanding, finding a single set of three cubes, we have found, using this approach, multiple sets of cubes. We have also given parametrisation for such cubes. Cases, where it is difficult or not possible to determine a seed equation, alternate method has been provided. Methods used are unattempted, innovative and easily comprehensible.

**Keywords:** Cube, Cubic Equation, Integer Root, Linear Equation, Parametrisation, Pythagorean's Triple, Quadratic Equation, Seed Equation.

**AMS Subject Classifications:** 11D25

## 1  Introduction

It is an open problem in *Number Theory* to characterise an integer as sum of three cubes. Heath-Brown conjectured that all integers except of the forms $9x \pm 4$ or $9x \pm 5$, can be expressed as sum of three cubes in multiple ways [4]. Apart from this statement, number 0 can not be expressed as sum of three distinct cubes on account of Fermat's last theorem that states, three integers $a, b$ and $c$ can not be expressed as $a^n + b^n = c^n$ where $n > 2$. Only possible solution is trivial, $a^3 + (-a)^3 + 0 = 0$. On the other hand, K. Mahler in the year 1936, found that integer 1 is expressible in infinite ways by parametric $(1 - 9p^3)^3 + (3p - 9p^4)^3 + (9p^4)^3 = 1$, [10]. There are other representations and family of parametrisation as well for this integer 1, [1]. Earlier in the year 1908, A.S. Verebrusov had already discovered that integer 2 is expressible in infinite ways by parametric $(1 - 6p^3)^3 + (1 + 6p^3)^3 + (-6p^2)^3 = 2$, [2],[8]. Notwithstanding parametric solutions stated above, there is a family of solutions for representing integer 2, [1], [4]. For integer 3, Mordell in the year 1953, expressed his opinion that he does not know anything more than to say $3 = 1^3 + 1^3 + 1^3 = 4^3 + 4^3 + (-5)^3$. Heath-Brown, Lioen, and te Riele in the year 1992, found the cubes that sum up to 39 by equation $39 = 134476^3 + 117367^3 + (-159380)^3$ [4]. Thereafter, for finding solutions to $x^3 + y^3 + z^3 = n$, where $n$ varies from 1 to 1000, many authors implemented computational searches up to $(|x|, |y|, |z|) < 10^{15}$ [5]. At present, up to $n < 100$, solutions have been found except two integers 33 and 42 [6]. But now, Andrew Booker and Andrew Southland have settled the cases for these two integers also [3], [11],[9]. Up to 1000, there are integers 114, 390, 627, 633, 732, 921 and 975 which could not be solved for their representation by sum of three cubes [1], [7].

Although much work has been done, we, resorting to innovative and unattempted methods, have endeavoured to solve the problem of representations of integers 1, 2, 3, ... so on except integers of the forms $9x \pm 4$ or $9x \pm 5$ by sum of three cubes. For integers of some specific forms, we have derived their parameterisation. We have also proved that a given integer can be expressed by multiple sets of three adding cubes. We have proved the veracity of formulae derived by us, presenting exhaustive examples in *Tables*. To begin with, we represent an integer $k$ as sum of three cubes by the equation
$$X^3 + Y^3 + Z^3 = k \tag{1.1}$$
where $X, Y$ and $Z$ are all integers positive or negative. To determine $X, Y$ and $Z$, we need a seed equation which, in fact, is a skeletal equation for integer $k$. This equation, then will be used to determine values of $X, Y$ and $Z$. Seed equations for $k = 1$, can be written easily as
$$(1)^3 + (0)^3 + (0)^3 = 1 \tag{1.2}$$
or
$$(1)^3 + (p)^3 + (-p)^3 = 1. \tag{1.3}$$
where $p$ is an integer positive or negative. Some of the seed equations for integers $1, 2, 3, ...$ are given in *Table* **1.1** for better clarity.

*Table* **1.1** Some Seed Equations

| Integers to Be Represented '$k$' | $A^3 + B^3 + C^3 = k$ | Integers to Be Represented '$k$' | $A^3 + B^3 + C^3 = k$ |
|---|---|---|---|
| 1 | $(1)^3 + (p)^3 + (-p)^3 = 1$ | 7 | $(2)^3 + (-1)^3 + (0)^3 = 7$ |
| 2 | $(1)^3 + (1)^3 + (0)^3 = 2$ | 8 | $(2)^3 + (0)^3 + (0)^3 = 8$ |
| 3 | $(1)^3 + (1)^3 + (1)^3 = 3$ | 9 | $(2)^3 + (1)^3 + (0)^3 = 9$ |
| 6 | $(2)^3 + (-1)^3 + (-1)^3 = 6$ | 10 | $(2)^3 + (1)^3 + (1)^3 = 10$ |

It may happen that seed equation for a given $k$ is not easily determinable, such cases will be dealt at the end. We submit, a real rational number, say $n$, can be expressed in algebraic form,
$$n = a \cdot x + A,$$
where $a$ and $A$ are real rational quantities assigned by us and $x$ is a variable. From above said equation, we can write, $x = \frac{n-A}{a}$. To illustrate, how a number, say $n$, can be expressed, we take the example for $n = 5$. Let $a$ and $A$ have values 7 and 3 respectively, then $x = \frac{5-3}{7} = \frac{2}{7}$. If $a$ and $A$ have values, 3 and $-1$ respectively, then $x = 2$. Integer 5 can, therefore, be written in infinite ways by assigning different values to $a$ and $A$. For example,
$$5 = 7 \cdot (2/7) + 3 = 3 \cdot 2 - 1 = \cdots.$$
That proves proceeding ***Lemma* 1.1.**

***Lemma* 1.1:** *A real rational number $n$ can always be expressed in the form $n = a \cdot x + A$, where $a$ and $A$ are real rational numbers and $x$ variable.*

### 1.1 Algebraic Equation for Sum of Cubes

Applying ***Lemma* 1.1,** $X, Y$ and $Z$ can be written in algebraic form, $X = a \cdot x + A$, $Y = b \cdot x + B$ and $Z = c \cdot x + C$. Equation (1.1), then can be written
$$(a \cdot x + A)^3 + (b \cdot x + B)^3 + (c \cdot x + C)^3 = k. \tag{1.4}$$
On expansion and rearrangement, this equation can be written,
$$x^3(a^3 + b^3 + c^3) + 3x^2(a^2 \cdot A + b^2 \cdot B + c^2 \cdot C) + 3x(a \cdot A^2 + b \cdot B^2 + c \cdot C^2)$$
$$+ A^3 + B^3 + C^3 - k = 0. \tag{1.5}$$
Since $a, b, c, A, B$ and $C$ have values as assigned by us, therefore, coefficients of $x^3, x^2, x$ and constant term of above equation are known. Being a cubic equation, finding integer values of $x$, satisfying this equation, is tedious, hence determination of values of $X = a \cdot x + A, Y = b \cdot x + B$ and $Z = c \cdot x + C$ will also be tedious and difficult. To tide over this difficulty, we will transform cubic into a linear equation.

### 2 Transformation of Cubic Equation into a Linear Equation

For solving cubic equation in an easy way, it is transformed into a quadratic equation and, then to linear equation. Transformation into a quadratic equation will be achieved by use of seed equation that has already been explained in the introduction section. Let the seed equation for integer $k$ be $A^3 + B^3 + C^3 = k$. Substituting this value of $k$ in Equation (1.5) and on simplifying, it gets transformed into
$$x\{x^2(a^3 + b^3 + c^3) + 3x(a^2 \cdot A + b^2 \cdot B + c^2 \cdot C) + 3(a \cdot A^2 + b \cdot B^2 + c \cdot C^2)\} = 0.$$
This equation has one root at $x = 0$ and this root is ignored on account of the fact, it will yield the given seed equation. Ignoring root at $x = 0$, the equation, then reduces to quadratic form
$$x^2(a^3 + b^3 + c^3) + 3x(a^2 \cdot A + b^2 \cdot B + c^2 \cdot C) + 3(a \cdot A^2 + b \cdot B^2 + c \cdot C^2) = 0. \quad (2.1)$$
If Equation (2.1) is solvable for real rational roots, say $x$ equal to $p$ and $q$ such that substitution of value of the roots in $(a \cdot x + A)$, $(b \cdot x + B)$ and $(c \cdot x + C)$ yields integer values on account of our choice of $a, b$ and $c$, then required cubes are given by relations
$$(a \cdot p + A)^3 + (b \cdot p + B)^3 + (c \cdot p + C)^3 = k, \quad (2.2)$$
and
$$(a \cdot q + A)^3 + (b \cdot q + B)^3 + (c \cdot q + C)^3 = k. \quad (2.3)$$
If the quadratic Equation (2.1) can not be solved or difficult to solve for real rational root or roots, it will be further transformed into a linear equation by equating constant term to zero. That is
$$(a \cdot A^2 + b \cdot B^2 + c \cdot C^2) = 0. \quad (2.4)$$
The transformed equation, then will be
$$x^2(a^3 + b^3 + c^3) + 3x(a^2 \cdot A + b^2 \cdot B + c^2 \cdot C) = 0.$$
Again this equation has one root, $x = 0$ and this root will yield seed equation, hence ignored. Transformed linear equation will, then be
$$x = -\frac{3(a^2 \cdot A + b^2 \cdot B + c^2 \cdot C)}{(a^3 + b^3 + c^3)} = r \text{ (say).} \quad (2.5)$$
Depending upon our selection of $a$, $b$ and $c$ if $r$ has a real rational value such that $(a \cdot x + A)$, $(b \cdot x + B)$ and $(c \cdot x + C)$ are integers, then solution will be
$$(a \cdot r + A)^3 + (b \cdot r + B)^3 + (c \cdot r + C)^3 = k. \quad (2.6)$$

## 2.1 Choice of $a$, $b$, $c$, $A$, $B$ and $C$

What values should be assigned to $A$, $B$ and $C$ depend upon the value of $k$ and that has already been explained while dealing with seed equations. Kindly refer to *Table* **1.1**. Assignment of values to $a$, $b$ and $c$ is very crucial. These are so chosen that either cubic equation (1.5) or quadratic equation (2.1) or linear equation (2.5) yields at least one real and rational root such that substitution of that root in $(a \cdot x + A)$, $(b \cdot x + B)$ and $(c \cdot x + C)$ gives integer values. To achieve this result, certain algebraic identities will be used. With this concept of transformation of cubic equation, we proceed to determine three cubes that sum up to a given integer.

## 3 Use of a Seed Equation

We will take up relatively easier cases first, then harder and ultimately the hardest, where it is difficult to find a seed equation.

### 3.1 Representation of Integer 2 and Other Integers Expressible as $2p^3$, Where $p$ Is an Integer

Seed equation for integer 2 is $(1)^3 + (1)^3 + (0)^3 = 2$. That makes $X = a \cdot x + 1$, $Y = b \cdot x + 1$ and $Z = c \cdot x$. We have chosen, $A = 1, B = 1$ and $C = 0$ intentionally so as to get rid of constant term in Equation (3.1). Equation, then can be written,
$$(a \cdot x + 1)^3 + (b \cdot x + 1)^3 + (c \cdot x)^3 = 2. \quad (3.1)$$
On simplification,
$$x^2(a^3 + b^3 + c^3) + 3x(a^2 + b^2) + 3(a + b) = 0.$$
For transformation into linear equation, we put $a + b = 0$. Transformed linear equation is
$$x(a^3 + b^3 + c^3) + 3(a^2 + b^2) = 0.$$
On putting, $b = -a$,
$$x = -\frac{6}{c}\left(\frac{a^2}{c^2}\right) = -\frac{6}{c}(p^2),$$

where $p = a/c$. Substituting this value of $x$ in Equation (2.6)

$$(1 - 6p^3)^3 + (1 + 6p^3)^3 + (-6p^2)^3 = 2. \qquad (3.2)$$

This parameterisation proves, there are infinite values of integer $p$ that satisfy Equation (3.2), hence there are infinite ways to represent 2 as sum of three cubes. This relation was discovered by A.S. Verebrusov who found this parametrisation [2]. But we have given the simplest method of determining family of sets of three cubes that sum up to 2. On multiplying identity (3.2) with $d^3$,

$$\{d(1 - 6p^3)\}^3 + \{d(1 + 6p^3)\}^3 + \{-d(6p^2)\}^3 = 2d^3. \qquad (3.3)$$

This identity is a parametrisation for the integers $2d^3$ expressed as sum of three cubes.

### 3.2 Representation of Integers 3 and $3p^3$

Seed equation for integer 3 is $(1)^3 + (1)^3 + (1)^3 = 3$. That makes $X = (a \cdot x + 1)$, $Y = (b \cdot x + 1)$ and $Z = (c \cdot x + 1)$. Corresponding equation is

$$(a \cdot x + 1)^3 + (b \cdot x + 1)^3 + (c \cdot x + 1)^3 = 3. \qquad (3.4)$$

On simplification,

$$x^3(a^3 + b^3 + c^3) + 3x^2(a^2 + b^2 + c^2) + 3x(a + b + c) = 0.$$

For transformation into linear equation, we put $a + b + c = 0$ and simplify it, then

$$x = 2\left(\frac{1}{c} - \frac{c}{ab}\right).$$

Substituting this value of $x$ in Equation (3.4), we get

$$\left\{1 + 2a\left(\frac{1}{c} - \frac{c}{ab}\right)\right\}^3 + \left\{1 + 2b\left(\frac{1}{c} - \frac{c}{ab}\right)\right\}^3 + \left\{1 + 2c\left(\frac{1}{c} - \frac{c}{ab}\right)\right\}^3 = 3. \qquad (3.5)$$

Let $b = c \cdot d$, then $a = -(c \cdot d + c)$ since $a = -(b + c)$, where $d$ is real rational quantity. On putting these values of $a$ and $b$ in Equation (3.5) and simplifying,

$$\left\{-1 - 2\left(d + \frac{1}{d}\right)\right\}^3 + \left\{1 + 2\left(d + \frac{1}{d+1}\right)\right\}^3 + \left\{1 + 2\left(\frac{1}{d} + \frac{d}{d+1}\right)\right\}^3 = 3. \qquad (3.6)$$

For integer solutions, $2\left(d + \frac{1}{d}\right)$, $2\left(d + \frac{1}{d+1}\right)$ and $2\left(\frac{1}{d} + \frac{d}{d+1}\right)$ must be integers. To achieve this objective, $d$ and $d + 1$ must be factors of 2 and either of the two, must not be zero. Due to reciprocity of $d$ and $d + 1$ with $1/d$ and $1/(d + 1)$ respectively as is apparent from Equation (3.6), values of $1/d$ and $1/(d + 1)$ will also satisfy Equation (3.6) if values of $d$ and $(d + 1)$ satisfy this equation. For example, $d = -2$ or $1/d = -2$ or $d = 1$ will satisfy the Equation (3.6) and yield

$$(4)^3 + (-5)^3 + (4)^3 = 3,$$
$$(4)^3 + (4)^3 + (-5)^3 = 3,$$
$$(-5)^3 + (4)^3 + (4)^3 = 3. \qquad (3.7)$$

These equations have same cubes but latter two equations have cubes displaced from their original position and, in this way, these three equations amount to one set of cubes. Multiplication of Equation (3.7) with $p^3$, where $p$ is an integer, gives parametrisation for $3p^3$

$$(4p)^3 + (4p)^3 + (-5p)^3 = 3p^3. \qquad (3.8)$$

### 3.2a Multiple Sets of Three Cubes That Represent Integer $3p^3$

Equation (3.6), when multiplied with $p^3$, where $p$ is an integer, takes the form

$$\left\{-p - 2p\left(d + \frac{1}{d}\right)\right\}^3 + \left\{p + 2p\left(d + \frac{1}{d+1}\right)\right\}^3 + \left\{p + 2p\left(\frac{1}{d} + \frac{d}{d+1}\right)\right\}^3 = 3p^3. \qquad (3.9)$$

For integer values of cubes, values of $d$ should be such that $2p(d + 1/d)$, $2p\{d + 1/(d + 1)\}$ and $2p\{1/d + d/(d + 1)\}$ must be integers. For integer 2, we have already found out factors in preceding paragraph and these will also be applicable to $2p$ as 2 is a factor of $2p$ but to have more factors, $2p$ must be of the form

$$2p = q(q + 1),$$

where $q$ is an integer positive or negative such that either $q$ or $(q + 1)$ is not zero. Let us take the case, where $q = 2$. In that case, factors of $2p$ (or 6) of the form $d(d + 1)$, are $(2)(3)$ and $(1)(2)$ ignoring factors $(-3)(-2)$ and $(-2)(-1)$ as these yield same cubes but displaced from their original positions. Therefore, putting $d = 1$ and $p = 3$, we get $(-15)^3 + (12)^3 + (12)^3 = 3(3)^3$ and on putting $d = 2$ and $p = 3$, we get $(-18)^3 + (17)^3 + (10)^3 = 3(3)^3$. This proves proceeding **Lemmas 3.1**, **3.2** and **3.3**.

**Lemma 3.1**: Integer 3 is represented by sum of cubes

$$\left\{-1 - 2\left(d + \frac{1}{d}\right)\right\}^3 + \left\{1 + 2\left(d + \frac{1}{d+1}\right)\right\}^3 + \left\{1 + 2\left(\frac{1}{d} + \frac{d}{d+1}\right)\right\}^3 = 3,$$

where $d$ has such values that $d$ and $d + 1$ are factors of $2$. Also $2\left(d + \frac{1}{d}\right)$, $2\left(d + \frac{1}{d+1}\right)$ and $2\left(\frac{1}{d} + \frac{d}{d+1}\right)$ have integer values.

**Lemma 3.2**: Integer $3p^3$ is represented by sum of cubes

$$\left\{-p - 2p\left(d + \frac{1}{d}\right)\right\}^3 + \left\{p + 2p\left(d + \frac{1}{d+1}\right)\right\}^3 + \left\{p + 2p\left(\frac{1}{d} + \frac{d}{d+1}\right)\right\}^3 = 3p^3,$$

where $d$ has such values that $d$ and $d + 1$ are factors of $2p$. Also $2p\left(d + \frac{1}{d}\right)$, $2p\left(d + \frac{1}{d+1}\right)$ and $2p\left(\frac{1}{d} + \frac{d}{d+1}\right)$ have integer values.

**Lemma 3.3**: If integer $2p$ has only factors $q$ and $(q + 1)$ and neither $q$ nor $(q + 1)$ is zero, then integer $3p^3$ is represented by two sets of sum of cubes as given below

$$\{-p - 2p(1 + 1)\}^3 + \left\{p + 2p\left(1 + \frac{1}{2}\right)\right\}^3 + \left\{p + 2p\left(1 + \frac{1}{2}\right)\right\}^3 = 3p^3$$

and

$$\left\{-p - 2p\left(q + \frac{1}{q}\right)\right\}^3 + \left\{p + 2p\left(q + \frac{1}{q+1}\right)\right\}^3 + \left\{p + 2p\left(\frac{1}{q} + \frac{q}{q+1}\right)\right\}^3 = 3p^3$$

provided $q$ is not further factorable as $q_1(q_1 + 1)$.

**3.2b Multiple Sets or Family of Sets of Three Cubes Representing Integer $3(P!/2!)^3$ for $P \geq 3$**

Let integer $p$ be such that $p = P!/2!$ and if $P = 3$, then $p = 3!/2! = 3$ and $2p = 6$. Six is factorable as $(1)(2)(3)$. Therefore, $(1)(2)$ and $(2)(3)$ are its two pairs of factors. That means $3(3!/2!)^3$ can be represented by two sets of three cubes.

Taking the case, when $P = 4$, then $p = 4!/2! = 12$ or $2p = 24$ which is factorable as $(1)(2)(3)(4)$, then $(1)(2)$, $(2)(3)$ and $(3)(4)$ are three pairs of factors of the form $d(d + 1)$. That means $3(4!/2!)^3$ can be represented by three sets of three cubes.

Takings another case, when $P = 5$ or $p = 5!/2! = 60$ or $2p = 120$ which is factorable as $(1)(2)(3)(4)(5)$, then $(1)(2)$, $(2)(3)$, $(3)(4)$, $(4)(5)$ and $(5)(6)$ are five pairs of factors of the form $d(d + 1)$. That means $3(5!/2!)^3$ can be represented by five sets of three cubes. Similarly, it can be proved, when $P = 6$ or $p = 6!/2! = 360$ or $2p = 720$, it is factorable by eight pairs of factors of the form $d(d + 1)$. That means $3(6!/2)^3$ can be represented by eight sets of three cubes by putting different values of $d$. These eight sets are given below.

$$(-1800)^3 + (1440)^3 + (1440)^3 = 3(6!/2!)^3 = 139968000$$
$$(-2160)^3 + (2040)^3 + (1200)^3 = 3(6!/2!)^3 = 139968000$$
$$(-2760)^3 + (2700)^3 + (1140)^3 = 3(6!/2!)^3 = 139968000$$
$$(-3420)^3 + (3384)^3 + (1116)^3 = 3(6!/2!)^3 = 139968000$$
$$(-4104)^3 + (4080)^3 + (1104)^3 = 3(6!/2!)^3 = 139968000$$
$$(-6210)^3 + (6200)^3 + (1090)^3 = 3(6!/2!)^3 = 139968000$$
$$(-6920)^3 + (6912)^3 + (1088)^3 = 3(6!/2!)^3 = 139968000$$
$$(-11208)^3 + (11205)^3 + (1083)^3 = 3(6!/2!)^3 = 139968000$$

From above, this can be concluded, when $k$ is of form form $3(P!/2!)^3$, where $P \geq 3$, then $k$ can be represented by multiple sets of three cubes. Number of such sets of three cubes equals to the number of factors of the form $d(d + 1)$ of $P!$. Also, when $P \geq 5$, then number of sets of cubes is more than $P$. This proves **Lemma 3.4**.

**Lemma 3.4:** *If integer $p = (P!/2!)$ and $P \geq 3$, then integer $k = 3(P!/2!)^3$ is represented by sum of cubes*

$$\left\{-p - 2p\left(d + \frac{1}{d}\right)\right\}^3 + \left\{p + 2p\left(d + \frac{1}{d+1}\right)\right\}^3 + \left\{p + 2p\left(\frac{1}{d} + \frac{d}{d+1}\right)\right\}^3 = 3p^3 = 3\left(\frac{P!}{2!}\right)^3$$

*in as many sets as $P!$ has distinct factors of the form $d(d + 1)$. When $P \geq 5$, then $3p^3$ can be represented by sets of cubes in as many numbers as equal to or more than $P$. Also, If $P!$ has $n$ sets of distinct factors, then using this method, there will be $n$ distinct sets of three cubes to represent $3p^3$.*

### 3.3 Representation of Integers 1 and Other Integers Expressible as $p^3$, Where $p$ Is an Integer

Seed equation for integer 1 is $(1)^3 + (0)^3 + (0)^3 = 1$. That makes $X = a \cdot x + 1$, $Y = b \cdot x$ and $Z = c \cdot x$. Therefore,

$$(a \cdot x + 1)^3 + (b \cdot x)^3 + (c \cdot x)^3 = 1. \tag{3.10}$$

On simplification,

$$x^2(a^3 + b^3 + c^3) + 3x(a^2) + 3(a) = 0.$$

Let $a = -9p^2$, $b = 3(1 - 3p^3)$, $c = 9p^3$, then quadratic equation gets transformed into

$$(1 - 9p^3)x^2 + 9xp^4 - p^2 = 0$$

which has roots, $x = p$ and $x = -p/(1 - p^3)$. On putting these values, $x = p$ or $x = -p/(1 - p^3)$ in Equation (3.10), we get the identity

$$(1 - 9p^3)^3 + (3p - 9p^4)^3 + (9p^4)^3 = 1. \tag{3.11}$$

This is a parameterisation for representation of integer 1 and proves, there are infinite values of $p$, hence infinite ways to represent 1 as sum of three cubes. This was first given by K. Mahler in 1936 [10]. But we have given the simplest method of determining family of sets of three cubes that sum up to 1. Referring to identity (3.11), if this is multiplied by $d^3$, where $d$ is an integer, it gets transformed into

$$\{d(1 - 9p^3)^3\}^3 + \{d(3 - 9p^3)\}^3 + \{d(9p^4)\}^3 = d^3 \tag{3.12}$$

which is a parametrisation for sum of three cubes that equals $d^3$.

### 3.4 Representation of Integer 4 or 5

Using **Lemma 1.1**, a real rational integer $X$ can be expressed as $X = 9x + a$, where $x$ is an integer positive or negative and $a$ is also an integer positive $0, 1, 2, 3, 4, 5$ or a negative integer $-1, -2, -3, -4$. For examination of divisibility of $X^3$ by 9, we expand $X^3$,

$$X^3 = (9x + a)^3 = (9x)^3 + 27x \cdot a(9x + a) + a^3.$$

In the above equation, $(9x)^3$ and $27x \cdot a(9x + a)$ are obviously divisible by 9, therefore, we will now discuss divisibility of $a^3$ by 9. Taking different values of $a$ as stated above, $a^3$ is found to have remainders either 0 or 1 or $-1$ or 0 or 1 or $-1$ or 0. Therefore, algebraic sum of three cubes i.e. $X^3 + Y^3 + Z^3$ will have a remainder which is either algebraic sum of any combination of three integers $0, 1$ or $-1$ or single integer 0 or 1 or $-1$ repeated thrice. All combinations of the sum of remainders yield 0 or $\pm 1$ or $\pm 2$ or $\pm 3$. Therefore, sum of three cubes when divided by 9 can not have remainders other than as mentioned above. Conversely, it can not have remainder $\pm 4$. When division is by 9, remainder $-4$ is equivalent to remainder 5 and remainder $-5$ equivalent to 4, therefore, it can also be stated that sum of three cubes, when divided by 9, can not have remainder $\pm 5$. Therefore, $k = 9x \pm 4$ or $k = 9x \pm 5$ can not be represented by sum of three cubes. In 1992, Roger Heath-Brown conjectured that every $k$ unequal to 4 or 5 modulo 9 has infinitely many representations as sums of three cubes [4].

### 3.5 Representation of Integer 7 and Other Integers Expressible as $7d^3$

Seed equation for integer 7 is $(2)^3 + (-1)^3 + (0)^3 = 7$. That makes $X = (a \cdot x + 2)$, $Y = (b \cdot x - 1)$, $Z = c \cdot x$ and

$$(a \cdot x + 2)^3 + (b \cdot x - 1)^3 + (c \cdot x)^3 = 7. \quad (3.13)$$

On putting $4a + b = 0$, $(c/a) = y$ and simplifying, transformed linear equation is

$$x = \frac{42a^2}{c^3 - 63a^3} = \frac{1}{a}\left(\frac{42}{y^3 - 63}\right) \quad (3.14)$$

On putting above said value of $x$, $b = -4a$ and $(c/a) = y$ in Equation (3.13), we get

$$\left\{2 + \frac{42}{y^3 - 63}\right\}^3 - \left\{1 + \frac{168}{y^3 - 63}\right\}^3 + \left\{\frac{42y}{y^3 - 63}\right\}^3 = 7.$$

At $y = 4$,

$$44^3 - 169^3 + 168^3 = 7.$$

On multiplying above said identity with $d^3$, it gets transforms into

$$(44d)^3 + (-169d)^3 + (168d)^3 = 7d^3$$

which is a parametrisation for integer expressible as $7d^3$, where $d$ is an integer.

### 3.6 Representation of an Integers Expressible as $p^3 - 1$ and $d^3(p^3 - 1)$ Where $d$ and $p$ Are Integers.

Seed equation for integer $p^3 - 1$ is $(p)^3 + (-1)^3 + (0)^3 = p^3 - 1$. That makes $X = (a \cdot x + p)$, $Y = (b \cdot x - 1)$, $Z = c \cdot x$ and

$$(a \cdot x + p)^3 + (b \cdot x - 1)^3 + (c \cdot x)^3 = p^3 - 1. \quad (3.15)$$

On expansion, simplification and putting $p^2 \cdot a + b = 0$, $(c/a) = y$, Equation 3.15 gets transformed into a linear equation,

$$x = \frac{3p \cdot a^2(p^3 - 1)}{c^3 - a^3(p^6 - 1)} = \left(\frac{1}{a}\right)\frac{3p(p^3 - 1)}{\frac{c^3}{a^3} - (p^6 - 1)} = \left(\frac{1}{a}\right)\frac{3p(p^3 - 1)}{y^3 - (p^6 - 1)} \quad (3.16)$$

On putting above said value of $x$, $b = -p^2 a$ and $(c/a) = y$ in Equation (3.15),

$$\left\{p + \frac{3p(p^3 - 1)}{y^3 - (p^6 - 1)}\right\}^3 + \left\{-1 - \frac{3p^3(p^3 - 1)}{y^3 - (p^6 - 1)}\right\}^3 + \left\{\frac{3p \cdot y(p^3 - 1)}{y^3 - (p^6 - 1)}\right\}^3 = p^3 - 1.$$

On putting, $y = p^2$ in above equation,

$$\{p + 3p(p^3 - 1)\}^3 + \{-1 - 3p^3(p^3 - 1)\}^3 + \{3p^3(p^3 - 1)\}^3 = p^3 - 1 \quad (3.17)$$

which is parametrisation for integer $p^3 - 1$ where $p$ is an integer positive or negative.

Table 3.1 Sum of Three Cubes That Equals to $p^3 - 1$

| $p$ | $X^3 + Y^3 + Z^3 = p^3 - 1$ | $p$ | $X^3 + Y^3 + Z^3 = p^3 - 1$ |
|---|---|---|---|
| 3 | $237^3 - 2107^3 + 2106^3 = 26$ | 8 | $12272^3 - 784897^3 + 784896^3 = 511$ |
| 4 | $760^3 - 12097^3 + 12096^3 = 63$ | 9 | $19665^3 - 1592137^3 + 1592136^3 = 728$ |
| 5 | $1865^3 - 46501^3 + 46500^3 = 124$ | 10 | $29980^3 - 2997001^3 + 2997000^3 = 999$ |
| 6 | $3876^3 - 139321^3 + 139320^3 = 215$ | 11 | $43901^3 - 5310691^3 + 5310690^3 = 1330$ |
| 7 | $7189^3 - 351919^3 + 351918^3 = 342$ | 12 | $62184^3 - 8952769^3 + 8952768^3 = 1727$ |

**Table 3.1** depicts sum of three cubes that equals to $p^3 - 1$ when integer $p$ varies from 3 to 12.

### 3.7 Representation of Integer 9, $9d^3$ and Other Integers Expressible as $p^3 + 1$

Examination of Equation (16) reveals that if $p$ is assumed as $-p$, this equation gets transformed into

$$\{p - 3p(p^3 + 1)\}^3 + \{1 + 3(p^3)(p^3 + 1)\}^3 + \{-3(p^3)(p^3 + 1)\}^3 = p^3 + 1. \quad (3.18)$$

At $p = 2$,

$$(-52)^3 + (217)^3 + (-216)^3 = 9.$$

On multiplying above said identity with $d^3$, it gets transformed into

$$(-52d)^3 + (217d)^3 + (-216d)^3 = 9d^3$$

which is a parametrisation for integers expressible as $9p^3$, where $p$ is an integer.

Table 3.2 Sum of Three Cubes That Equals to $p^3 + 1$

| p | $X^3 + Y^3 + Z^3 = p^3 + 1$ | p | $X^3 + Y^3 + Z^3 = p^3 + 1$ |
|---|---|---|---|
| 3 | $-249^3 + 2269^3 - 2268^3 = 28$ | 8 | $-12304^3 + 787969^3 - 787968^3 = 513$ |
| 4 | $-776^3 + 12481^3 - 12480^3 = 65$ | 9 | $-19701^3 + 1596511^3 - 1596510^3 = 730$ |
| 5 | $-1885^3 + 47251^3 - 47250^3 = 126$ | 10 | $-30020^3 + 3003001^3 - 3003000^3 = 1001$ |
| 6 | $-3900^3 + 140617^3 - 140616^3 = 217$ | 11 | $-43945^3 + 5311867^3 - 5318676^3 = 1332$ |
| 7 | $-7217^3 + 353977^3 - 353976^3 = 344$ | 12 | $-62232^3 + 8963137^3 - 8963136^3 = 1729$ |

Table 3.2 depicts sum of three cubes that equals to $p^3 + 1$ where $p$ is an integer from 3 to 12.

### 3.8 Representation of Integer 10 and Other Integers Expressible as $10p^3$

Seed equation for integer 10 is $(2)^3 + (1)^3 + (1)^3 = 10$. That makes $X = (a \cdot x + 2)$, $Y = (b \cdot x + 1)$, $Z = (c \cdot x + 1)$ and

$$(a \cdot x + 2)^3 + (b \cdot x + 1)^3 + (c \cdot x + 1)^3 = 10. \quad (3.19)$$

On putting, $4a + b + c = 0$, $b = c = -2a$ and simplifying,

$$x = 2/a \quad (3.20)$$

and on putting this value of $x$ in Equation (3.19), we get

$$(4)^3 + (-3)^3 + (-3)^3 = 10.$$

Therefore, parametrisation for $10p^3$ is

$$(4p)^3 + (-3p)^3 + (-3p)^3 = 10p^3.$$

### 3.9 Representation of Integers Using Pythagorean's Triples

Consider $A, B$ and $C$ as Pythagorean's Triple such that $A^2 + B^2 = C^2$ and its seed equation as $k = A^3 - B^3 + C^3$. Therefore, $X = (a \cdot x + A)$, $Y = (a \cdot x - B)$ and $Z = (-a \cdot x + C)$ and

$$(a \cdot x + A)^3 + (a \cdot x - B)^3 + (-a \cdot x + C)^3 = A^3 - B^3 + C^3. \quad (3.21)$$

On simplification by putting $A^2 + B^2 - C^2 = 0$, being Pythagorean's triple,

$$x = -(3/a)(A - B + C).$$

Substituting this value of $x$ in Equation (3.21)

$$\{-3(A - B + C) + A\}^3 + \{-3(A - B + C) - B\}^3 + \{3(A - B + C) + C\}^3 = A^3 - B^3 + C^3.$$

Let $S = B - (A + C)$, then the Equation takes the form,

$$(3S + A)^3 + (3S - B)^3 + (-3S + C)^3 = A^3 - B^3 + C^3. \quad (3.22)$$

Since $(y^2 - 1)$, $(2y)$ and $(y^2 + 1)$ are also Pythagorean's triple for all integer values of $y$, therefore, substituting these for $A, B$ and $C$ respectively, identity (3.22) takes the form,

$$(-5y^2 + 6y - 1)^3 + (-6y^2 + 4y)^3 + (7y^2 - 6y + 1)^3 = 2y^2(y^4 - 4y + 3). \quad (3.23)$$

At $y = -1, 2$ and 3, we obtain following sum of cubes

$$(-12)^3 + (-10)^3 + (14)^3 = 2(8),$$
$$(-9)^3 + (-16)^3 + (17)^3 = 8(11),$$
$$(-28)^3 + (-42)^3 + (46)^3 = 18(72).$$

respectively. On putting, $p = 1/y$ in Equation (3.23),

$$(-p^2 + 6p - 5)^3 + (4p - 6)^3 + (p^2 - 6p + 7) = 2(3p^4 - 4p^3 + 1). \quad (3.24)$$

At $p = 2, 7$ and 4, we obtain following sum of cubes

$$(3)^3 + (2)^3 + (-1)^3 = 34,$$
$$(-12)^3 + (22)^3 + (14)^3 = 11664,$$
$$(3)^3 + (10)^3 + (-1)^3 = 1026.$$

respectively. Also, since $A = p^2 - q^2$, $B = 2p \cdot q$, $C = p^2 + q^2$ are Pythagorean's triple for all real rational values of $p$ and $q$, therefore,

$$\{-(5p - q)(p - q)\}^3 + \{2p(2q - 3p)\}^3 + \{7p^2 + q^2 - 6p \cdot q\}^3$$

$$= 2p^2(p^4 + 3q^4 - 4p \cdot q^3). \tag{3.25}$$

At $p, q = (1, 0)$, $p, q = (1, 2)$ and $p, q = (2, 1)$, we obtain following sum of cubes

$$(-5)^3 + (-6)^3 + (7)^3 = 2,$$
$$(3)^3 + (2)^3 + (-1)^3 = 2(17),$$
$$(-9)^3 + (-16)^3 + (17)^3 = 8(11).$$

respectively. ***Table 3.3*** depicts sum of cubes for $k = A^3 - B^3 + C^3$ where $A, B$ and $C$ are Pythagorean's triple

***Table 3.3*** Sum of Cubes for $k = A^3 - B^3 + C^3$

| $A^3 - B^3 + C^3 = k$ | $(3S + A)^3 + (3S - B)^3 - (3S - C)^3$ | $A^3 - B^3 + C^3 = k$ | $(3S + A)^3 + (3S - B)^3 - (3S - C)^3$ |
|---|---|---|---|
| $1^3 - 0^3 + 1^3 = 2$ | $-5^3 - 6^3 + 7^3$ | $5^3 - 12^3 - 13^3 = -3800$ | $65^3 + 48^3 - 73^3$ |
| $3^3 - 4^3 - 5^3 = -162$ | $21^3 + 14^3 - 23^3$ | $5^3 + 12^3 + 13^3 = 4050$ | $-85^3 - 78^3 + 103^3$ |
| $3^3 - 4^3 + 5^3 = 88$ | $-9^3 - 16^3 + 17^3$ | $8^3 + 15^3 - 17^3 = -1026$ | $-10^3 - 3^3 + 1^3$ |
| $3^3 + 4^3 + 5^3 = 216$ | $-33^3 - 32^3 + 41^3$ | $8^3 + 15^3 + 17^3 = 8800$ | $-112^3 - 105^3 + 137^3$ |
| $3^3 + 4^3 - 5^3 = -34$ | $-3^3 - 2^3 + 1^3$ | $8^3 - 15^3 + 17^3 = 2050$ | $-22^3 - 45^3 + 47^3$ |
| $5^3 - 12^3 + 13^3 = 594$ | $-13^3 - 30^3 + 31^3$ | $8^3 - 15^3 - 17^3 = -7776$ | $80^3 + 57^3 - 89^3$ |

## 4 Parametrisation for Integers Conforming to Some Specific Forms

### 4.1 Parametrisation When Given Integer Is of Form $k = 9p^3 + 1$

When $k = 9p^3 + 1$, we will consider seed equation $(1)^3 + (0)^3 + (0)^3 = 1$ and for $9p^3$, we will find this value from coefficient of $x$ after proper assignment of values to $a, b$ and $c$. For that, consider $X = (a \cdot x + 1)$, $Y = (-a \cdot x)$, $Z = (c \cdot x)$ and $X^3 + Y^3 + Z^3 = (d \cdot x + 1)$, then

$$(a \cdot x + 1)^3 + (-a \cdot x)^3 + (c \cdot x)^3 = d \cdot x + 1. \tag{4.1}$$

On putting, $d = 3a$, $a/c = p$ and simplifying, it takes the form

$$3ax(ax) + c^3 x^3 = 0,$$
$$x = -\frac{3}{c}(p^2).$$

Substituting this value of $x$ in Equation (4.1),

$$\{3p^3 - 1\}^3 + \{-3p^3\}^3 + \{3p^2\}^3 = 9p^3 - 1. \tag{4.2}$$

At $p = 2$, 3 and 4, following are sum of cubes

$$23^3 - 24^3 + 12^3 = 9(2^3) - 1 = 71,$$
$$80^3 - 81^3 + 27^3 = 9(3^3) - 1 = 242,$$
$$191^3 - 192^3 + 48^3 = 9(4^3) - 1 = 575.$$

Substituting $-p$ for $p$

$$\{(3p^3 + 1)\}^3 + \{-3p^3\}^3 + \{-3p^2\}^3 = 9p^3 + 1. \tag{4.3}$$

At $p = 1$, 2 and 3, following are sum of cubes

$$4^3 - 3^3 - 3^3 = 9(1^3) + 1 = 10,$$
$$25^3 - 24^3 - 12^3 = 9(2^3) + 1 = 73,$$
$$82^3 - 81^3 - 27^3 = 9(4^3) - 1 = 244.$$

### 4.2 Parametrisation, When Given Integer Is of Form $k = 9(p + 1)(p^2 - 1)$

We will consider seed equation for integer 0 as $(1)^3 + (-1)^3 + (0)^3 = 0$ and for $k = 9(p + 1)(p^2 - 1)$, we will assign values to $a, b$ and $c$ so as to obtain this value of $k$ from coefficient of $x$ while transforming it to linear equation. Consider $X = (a \cdot x + 1)$, $Y = (c \cdot x - 1)$, $Z = (-a \cdot x)$ and $X^3 + Y^3 + Z^3 = dx$, then

$$(a \cdot x + 1)^3 + (c \cdot x - 1)^3 - (a \cdot x)^3 = d \cdot x. \tag{4.4}$$

On expanding and rearranging, it takes the form
$$c^3 \cdot x^2 + 3x^2(a^2 - c^2) + x(3a + 3c - d) = 0.$$

On putting, $d = 3(a + c)$ and $\frac{a}{c} = p$ for reducing it to linear equation, it takes the form
$$x = -\frac{3}{c}(p^2 - 1).$$

On putting this value of $x$ in Equation (4.4), it takes the form
$$\{3p(p^2 - 1) - 1\}^3 + \{3(p^2 - 1) + 1\}^3 + \{-3p(p^2 - 1)\}^3 = 9(p + 1)(p^2 - 1). \tag{4.5}$$

At $p = 2$, 3 and 4, we obtain
$$(17)^3 + (10)^3 + (-18)^3 = 81,$$
$$71^3 + 25^3 - 72^3 = 9(4)(8) = 288,$$
$$179^3 + 46^3 - 180^3 = 9(5)(15) = 675.$$

Substituting $-p$ for $p$ and rearranging, Equation (4.5) is into
$$\{3p(p^2 - 1) + 1\}^3 - \{3(p^2 - 1) + 1\}^3 - \{3p(p^2 - 1)\}^3 = 9(p - 1)(p^2 - 1). \tag{4.6}$$

At $p = 2$, 3 and 4, we obtain
$$19^3 - 10^3 - 18^3 = 9(1)(3) = 27,$$
$$(73)^3 - (25)^3 + (-72)^3 = 9(2)(8) = 144,$$
$$181^3 - 46^3 - 180^3 = 9(3)(15) = 405.$$

### 4.3 Parametrisation, When Given Integer Is of Form $k = (9p^3 - 1)^2$ and $(9p^3 + 1)^2$

Consider $X = (a \cdot x + 1)$, $Y = (-a \cdot x)$, $Z = (c \cdot x)$ and $X^3 + Y^3 + Z^3 = (d \cdot x + 1)^2$, then
$$(a \cdot x + 1)^3 - (a \cdot x)^3 + (c \cdot x)^3 = (d \cdot x + 1)^2. \tag{4.7}$$

On expansion and simplification,
$$x^3 \cdot c^3 + 3x^2(a) + 3x \cdot a^2 = d^2 \cdot x^2 + 2d \cdot x.$$

For transforming it into a linear equation, we put $d = 2$ and $a = 4/3$, then $x = -4/(3c^3)$. On substituting this value of $x$ in Equation (4.7) and putting $p = 2/3c$,
$$(-6p^3 + 1)^3 + (6p^3)^3 + (-3p^2)^3 = (9p^3 - 1)^2. \tag{4.8}$$

At $p = 1$, 2 and 3, we obtain
$$(-5)^3 + (6)^3 + (-3)^3 = (9 - 1)^2 = 64,$$
$$(-47)^3 + (48)^3 + (-12)^3 = (71)^2 = 5041,$$
$$(-161)^3 + (162)^3 + (-27)^3 = (242)^2 = 58564.$$

Substituting $-p$ for $p$ and rearranging
$$(6p^3 + 1)^3 + (-6p^3)^3 + (-3p^2)^3 = (9p^3 + 1)^2. \tag{4.9}$$

At $p = 1$, 2 and 3, we obtain
$$(7)^3 + (-6)^3 + (-3)^3 = (10)^2 = 100,$$
$$(49)^3 + (-48)^3 + (-12)^3 = (73)^2 = 5329,$$
$$(163)^3 + (-162)^3 + (-27)^3 = (244)^2 = 59536.$$

### 4.4 Parametrisation, When Given Integer Is of Form $k = 2(5 - 9p^2)$

Consider $X = (a \cdot x + 1)$, $Y = (-a \cdot x + 1)$, $Z = (-c \cdot x - 1)$ and $X^3 + Y^3 + Z^3 = d \cdot x + 1$, then
$$(a \cdot x + 1)^3 + (-a \cdot x + 1)^3 + (-c \cdot x - 1)^3 = d \cdot x + 1. \tag{4.10}$$

On expanding, then putting $d = -3c$ and $p = a^2/c^2$, it gets transformed into linear equation
$$x = 3\left(\frac{2a^2 - c^2}{c^3}\right) = \frac{3}{c}\left(\frac{2a^2}{c^2} - 1\right) = \frac{3}{c}(2p^2 - 1).$$

Substitution of this value of $x$ in Equation (4.10) yields
$$\{6p^3 - 3p + 1\}^3 + \{-6p^3 + 3p + 1\}^3 + \{-6p^2 + 2\}^3 = 2\{5 - 9p^2\}. \tag{4.11}$$

At $p = 2$, 3 and 4, we obtain
$$(43)^3 + (-41)^3 + (-22)^3 = (-72 + 10) = -62,$$

$$(154)^3 + (-152)^3 + (-52)^3 = (-162 + 10) = -152,$$
$$(373)^3 + (-371)^3 + (-94)^3 = (-288 + 10) = -278.$$

### 4.5 Miscellaneous Parametrisation

Consider the equation
$$(a \cdot x + 1)^3 + (-a \cdot x + 1)^3 + (-c \cdot x - C)^3 = (d \cdot x + 2 - C^3). \quad (4.12)$$

On expansion, simplification, and putting $d = -3c \cdot C^2$, it gets transformed into linear equation $x = 3\left(\frac{2a^2 - C \cdot c^2}{c^3}\right)$. Substituting this value of $x$ and putting, $\frac{a}{c} = p$ in Equation (4.12),

$$\{3(2p^3 - C \cdot p) + 1\}^3 - \{3(2p^3 - C \cdot p) - 1\}^3 - \left\{\frac{3}{p}(2p^3 - C \cdot p) + C\right\}^3$$
$$= \left\{-\left(\frac{9C^2}{p}\right)(2p^3 - Cp) + 2 - C^3\right\}. \quad (4.13)$$

When $C = p^2$, Equation (4.13) gets transformed into
$$\{3p^3 + 1\}^3 - \{3p^3 - 1\}^3 - \{4p^2\}^3 = \{-10p^6 + 2\}. \quad (4.14)$$

At $p = 2, 3$ and $4$, we obtain
$$(25)^3 + (-23)^3 + (-16)^3 = -640 + 2 = -638,$$
$$(82)^3 + (-80)^3 + (-36)^3 = -7290 + 2 = -7288,$$
$$(193)^3 + (-191)^3 + (-64)^3 = -40960 + 2 = -40958.$$

When $C = p$, Equation (4.13) gets transformed into
$$\{3p^2(2p - 1) + 1\}^3 - \{3p^2(2p - 1) - 1\}^3 - \{2p(3p - 1)\}^3 = 2(1 + 4p^3 - 9p^4). \quad (4.15)$$

At $p = -1, 2$ and $-2$, we obtain
$$(-8)^3 + (10)^3 + (-8)^3 = 2(-12) = -24,$$
$$(37)^3 + (-35)^3 + (-20)^3 = 2(-111) = -222,$$
$$(-59)^3 + (61)^3 + (-28)^3 = 2(-175) = -350.$$

When $C = -1$, Equation (4.13) gets transformed into
$$\{1 + 3p(2p^2 + 1)\}^3 - \{3p(2p^2 + 1) - 1\}^3 - \{3(2p^2 + 1) - 1\}^3 = -\{6 + 18p^2\}. \quad (4.16)$$

At $p = 0, 2$ and $3$, we obtain
$$(1)^3 + (1)^3 + (-2)^3 = -6,$$
$$(55)^3 + (-53)^3 + (-26)^3 = -78,$$
$$(172)^3 - (170)^3 - (56)^3 = -168,$$

When $C = -2$, Equation (4.13) gets transformed into
$$\{1 + 6p(p^2 + 1)\}^3 - \{6p(p^2 + 1) - 1\}^3 - \{6(p^2 + 1) - 2\}^3 = -\{62 + 72p^2\}. \quad (4.17)$$

At $p = 0, 1$ and $2$, we obtain
$$(1)^3 + (1)^3 + (-4)^3 = -62,$$
$$(13)^3 + (-11)^3 + (-10)^3 = -(62 + 72) = -134,$$
$$(61)^3 + (-59)^3 + (-28)^3 = -(62 + 288) = -350.$$

When $C = 2$, Equation (4.13) gets transformed into
$$\{1 + 6p(p^2 - 1)\}^3 - \{6p(p^2 - 1) - 1\}^3 - \{6(p^2 - 1) + 2\}^3 = 6\{11 - 12p^2\}, \quad (4.18)$$

At $p = 2, 3$ and $4$, we obtain
$$(37)^3 + (-35)^3 + (-20)^3 = 6(11 - 48) = -222,$$
$$(145)^3 + (-143)^3 + (-50)^3 = 6(11 - 108) = -582,$$
$$(361)^3 + (-359)^3 + (-92)^3 = 6(11 - 192) = -1086.$$

When $C = 1/p$, Equation (4.13) gets transformed into
$$\{p(3p^3 - 1)\}^3 - \{p(3p^3 - 2)\}^3 - \{(3p^3 - 1)\}^3 = 1 - 2p^3. \quad (4.19)$$

At $p = -1, 2$ and $-2$, we obtain
$$(4)^3 + (-5)^3 + (4)^3 = 3,$$
$$(46)^3 + (-44)^3 + (-23)^3 = -15,$$

$$(50)^3 + (-52)^3 + (25)^3 = 17.$$

When $C = 1/2p$, Equation (4.13) gets transformed into
$$\{p(12p^3 - 1)\}^3 - \{p(12p^3 - 5)\}^3 - \{2(6p^3 - 1)\}^3 = 4(2 - 5p^3). \tag{4.20}$$

At $p = 1, -1$ and 2, we obtain
$$(11)^3 + (-7)^3 + (-10)^3 = -12,$$
$$(13)^3 + (-17)^3 + (14)^3 = 28,$$
$$(190)^3 + (-182)^3 + (-94)^3 = -152.$$

On putting $p = 2t$ in Equation (4.13), we obtain
$$\{t(96t^3 - 1)\}^3 - \{t(96t^3 - 5)\}^3 - \{(48t^3 - 1)\}^3 = (1 - 20t^3). \tag{4.21}$$

At $t = 1, -1$ and 2,
$$95^3 - 91^3 - 47^3 = -19,$$
$$97^3 - 101^3 + 49^3 = 21,$$
$$1534^3 - 1526^3 - 383^3 = -159.$$

In general, consider equation
$$(a \cdot x + q)^3 - (a \cdot x - q)^3 - (c \cdot x + C)^3 = (d \cdot x + 2q^3 - C^3). \tag{4.22}$$

On expanding, simplifying, putting, $d = -3c \cdot C^2$ and $\frac{a}{c} = p$, it gets transformed into linear equation
$$x = \frac{3}{c}(2q \cdot p^2 - C).$$

Substituting this value of $x$ in Equation (4.22),
$$\{3p(2q \cdot p^2 - C) + q\}^3 - \{3p(2q \cdot p^2 - C) - q\}^3 - \{3(2q \cdot p^2 - C) + C\}^3$$
$$= 2C^2(4C - 9p^2 \cdot q) + 2q^3. \tag{4.23}$$

We consider one more general equation
$$(a \cdot x + r)^3 - (a \cdot x - q)^3 - (c \cdot x + C)^3 = (d \cdot x + q^3 + r^3 - C^3). \tag{4.24}$$

On expanding, simplifying, putting, $\frac{a}{c} = p$ and $d = 3c\{p(r^2 - q^2) - C^2\}$, it gets transformed into
$$x = \frac{3}{c}\{p^2(r + q) - C\}.$$

Substituting this value of $x$ in Equation (4.24),
$$[3p\{p^2(r + q) - C\} + r]^3 - [3p\{p^2(r + q) - C\} - q]^3 - [\{3p^2(r + q) - 2C\}]^3$$
$$= 9\{p^2(r + q) - C\}\{p(r^2 - q^2) - C^2\} + r^3 + q^3 - C^3. \tag{4.25}$$

On putting different values of $p, r, q$ and $C$, we will get parametrisation for various values of $k$.

### 4.6 Parametric Solutions to Representation of Integers of the Form ±2 Modulo 9 and ±3 Modulo 9 by Sum of Three Cubes

We have already taken up the cases, when integer $k = 9x \pm 1$ in paragraphs 4.1 to 4.6. Although a few cases of integer $k = 9x \pm 2$ and $k = 9x \pm 3$ have also been considered in aforementioned paragraphs, we, in order, to further elaborate, take up these cases in this sub paragraph.

a. When Integers are of the form ±2 modulo 9.

Let $(a \cdot x + 1)^3 + (-a \cdot x)^3 + (c \cdot x + 1)^3 = d \cdot x + 2$. On simplifying, putting $d = 3(a + c)$ and $p = a/c$, we get, $x = -(3/c)(1 + p^2)$. On substituting this value of $x$ in the equation, $(a \cdot x + 1)^3 + (-a \cdot x)^3 + (c \cdot x + 1)^3 = d \cdot x + 2$ and, then rearranging, we get parametric solution,
$$\{-3p(1 + p^2) + 1\}^3 + \{3p(1 + p^2)\}^3 + \{-3(1 + p^2) + 1\}^3 = -9(p + 1)(1 + p^2) + 2$$

Let $p = -q$, then
$$\{3q(1 + q^2) + 1\}^3 + \{-3q(1 + q^2)\}^3 + \{-3(1 + q^2) + 1\}^3 = 9(q - 1)(1 + q^2) + 2. \tag{4.26}$$

At $q = 2, 3$ and 4, we get
$$(31)^3 + (-30)^3 - (14)^3 = 47,$$
$$(91)^3 + (-90)^3 + (-29)^3 = 182,$$
$$(205)^3 + (-204)^3 + (-50)^3 = 461$$

respectively. In addition to the above parametric solution, we can have another parametric solution by putting $r = 7, q = -5$ and $C = 6$ in Equation (4.25) which gets transformed into
$$[3p\{2p^2 - 6\} + 7]^3 - [3p\{2p^2 - 6\} + 5]^3 - [\{6p^2 - 12\}]^3 = 216\{p^2 - 3\}\{2p - 3\} + 2. \quad (4.27)$$
At $p = 2, 1, 3$, we get
$$(19)^3 + (-17)^3 + (-12)^3 = 218,$$
$$(-5)^3 + (7)^3 + (6)^3 = 434$$
$$(115)^3 + (-113)^3 + (42)^3 = 3890$$
respectively. Thus Equations (4.26) and (4.27) yield parametric solution to represent integer $k = 9x + 2$ by sum of three cubes.

Assuming $p$ to be $-p$, Equation (4.27) gets transformed into
$$[3p\{2p^2 - 6\} - 7]^3 - [3p\{2p^2 - 6\} - 5]^3 + [\{6p^2 - 12\}]^3 = 216\{p^2 - 3\}\{2p + 3\} - 2. \quad (4.28)$$
At $p = 2, 3, 4$, we get
$$(5)^3 + (-7)^3 + (12)^3 = 1510,$$
$$(101)^3 + (-103)^3 + (42)^3 = 11662,$$
$$(305)^3 + (-307)^3 + (84)^3 = 30884$$
respectively. Thus Equations (4.28) yields parametric solution to represent integer $k = 9x - 2$ by sum of three cubes.

b. *When Integers are of the form $\pm 3$ modulo 9.*

Using Equation (4.25) and putting $r = 1, q = 1, C = -1$, we get
$$-[3p\{2p^2 + 1\} + 1]^3 + [3p\{2p^2 + 1\} - 1]^3 + [6p^2 + 2]^3 = 9\{2p^2 + 1\} - 3. \quad (4.29)$$
At $p = 1, 2, 3$, we get
$$(-10)^3 + (8)^3 + (8)^3 = 24,$$
$$(-55)^3 + (53)^3 + (26)^3 = 78,$$
$$(-172)^3 + (170)^3 + (56)^3 = 168$$
respectively. Thus Equations (4.29) yields parametric solution to represent integer $k = 9x - 3$ by sum of three cubes.

Substituting $p$ with $-(3q + 1)$ in Equation (4.19), we get
$$[3(3q + 1)^4 + 3q + 1]^3 - [3(3q + 1)^4 + 6q + 2]^3 + [3(3q + 1)^3 + 1]^3 = 2(3q + 1)^3 + 1. \quad (4.30)$$
This is a parametrisation for representing integer of the form $9x + 3$. At $q = 0, 1, 2$, we get
$$(4)^3 + (-5)^3 + (4)^3 = 3,$$
$$(772)^3 + (-776)^3 + (193)^3 = 129,$$
$$(7210)^3 + (-7217)^3 + (1030)^3 = 687.$$
respectively.

## 5 Method When Seed Equation Is Not Helpful

There may be cases, where use of seed equation may not be fruitful in determining cubes by transformation of cubic into linear equation. For example, it is per se difficult to find three cubes that sum up to 6 except its seed equation $(2)^3 + (-1)^3 + (-1)^3 = 6$. For finding another set of three cubes, we consider $X = (a \cdot x + 2)$, $Y = (b \cdot x - 1)$, $Z = c \cdot x - 1$, therefore,
$$(a \cdot x + 2)^3 + (b \cdot x - 1)^3 + (c \cdot x - 1)^3 = 6. \quad (5.1)$$
On expansion and simplification, it gets transforms into quadratic equation
$$x^2(a^3 + b^3 + c^3) + 3x(2a^2 - b^2 - c^2) + 3(4a + b + c) = 0. \quad (5.2)$$
On transforming it into linear equation, we put $b + c = -4a$ and $c/a = y$. Resultant equation is
$$x = -\left(\frac{2}{a}\right)\frac{y^2 + 4y + 7}{4y^2 + 16y + 21} = -\left(\frac{2}{a}\right)z.$$
where $z = \frac{y^2+4y+7}{4y^2+16y+21}$. On putting this value of $z$ in Equation (5.1), we get
$$\{-4z + 2\}^3 + \{2z(4 + y) - 1\}^3 + \{-2y \cdot z - 1\}^3 = 6. \quad (5.3)$$

For ensuring adding cubes to be integers, value of $y$ should be such that $z$ must have either integer value or a value of the form $p/2$, where $p$ is an integer. This value is not easily found out. To determine this value of $y$, different integer values are assigned one by one to $y$ and at each value, it is checked whether resultant value of $z$ satisfies our requirement. Consider $p_1$ as that value of $z$ say at $y = y_1$, then
$$\{-4p_1 + 2\}^3 + \{2p_1(4 + y_1) - 1\}^3 + \{-2y_1 \cdot z_1 - 1\}^3 = 6.$$
Alternatively, in stead of integer values, we can either assign values to $y$ as $y_1/2$ or $y_1$ where $y_1$ is an integer so that $z$ must, then have integer value. Purpose of such assignments is that $\{-4z + 2\}$, $\{2z(4 + y) - 1\}$ and $\{-2y \cdot z - 1\}$ must have integer values.

In above method, we have put a condition that $b + c = -4a$ so that quadratic (5.2) transforms into a linear equation. However, a general solution can also be found out by solving the quadratic equation for its real roots. Let $a = a_1 \cdot c$ and $b = b_1 \cdot c$ where $a_1$ and $b_1$ are real rational quantities, then quadratic equation takes the form
$$c^2 \cdot x^2(a_1^3 + b_1^3 + 1) + 3c \cdot x(2a_1^2 - b_1^2 - 1) + 3(4a_1 + b_1 + 1) = 0 \quad (5.4)$$
which has roots
$$c \cdot x = \frac{-Q \pm (Q^2 - 4P \cdot R)^{\frac{1}{2}}}{2P} = z \text{ (say)}, \quad (5.5)$$
where $P = (a_1^3 + b_1^3 + 1)$, $Q = 3(2a_1^2 - b_1^2 - 1)$ and $R = 3(4a_1 + b_1 + 1)$. Substituting the value of $x$ given by Equation (5.5) in Equation (5.1), we get
$$(a_1 \cdot z + 2)^3 + (b_1 \cdot z - 1)^3 + (z - 1)^3 = 6.$$
That requires $a_1$ and $b_1$ should have such real rational values that $z$, $(a_1 \cdot z)$ and $(b_1 \cdot z)$ given by Equation (5.5) must have integer values. Finding such values of $a_1$ and $b_1$ requires assumption of different values for these till the condition is satisfied.

## 6 Method When Seed Equation Is Not Determinable

So far, we dealt with those cases where seed equation was easily determinable. But there are cases, where it is not possible to find or difficult to find a seed equation for a given integer $k$. In such cases, an equation will be found for arbitrary integer $k_1$ which is different from $k$. Let the equation for $k_1$ is $a_1^3 + b_1^3 + c_1^3 = k_1$. Please mind this equation unlike the seed equations already described, will not transform cubic equation into quadratic equation. Sum of cubes can, then be written as
$$(a \cdot x + a_1)^3 + (b \cdot x + b_1)^3 + (c \cdot x + c_1)^3 = k. \quad (6.1)$$
On expansion and simplification by assuming $a = a_2 \cdot c$, $b = b_2 \cdot c$ and $c \cdot x = z$, where $a_2$ and $b_2$ are real rational quantities,
$$z^3(a_2^3 + b_2^3 + 1) + 3z^2(a_1 \cdot a_2^2 + b_1 \cdot b_2^2 + c_1) + 3z(a_1^2 \cdot a_2 + b_1^2 \cdot b_2 + c_1^2) - d = 0 \quad (6.2)$$
where $k - k_1 = d$. Equation (6.2) is a cubic equation which can have roots either (1) two complex and one real or (2) two irrational and one rational root or (3) three real rational roots. Nature of roots depends upon values of coefficients of $z^3$, $z^2$, $z$ and constant term. Assignment of different integer values to $a_2$, $b_2$ and $c$ will give different roots and the root (or roots) that has (or have) integer value(s) will be the required root(s). Let that integer root be $c \cdot x = z$ or $x = z/c$, then required cubes will be
$$(a_2 \cdot z + a_1)^3 + (b_2 \cdot z + b_1)^3 + (z + c_1)^3 = k \quad (6.3)$$
Different procedures and algorithms have been devised by mathematicians to find the desired roots but we are not going into details of these as the purpose of this paper is to represent a given integer by sum of three cubes using novel approach of seed equation.

## 7 Multiple Representations as Sums of Three Cubes for Integer, Where $k \neq 4$ or $k \neq 5$ Modulo 9

It is reiterated, method of seed equation employs use of a seed or skeletal equation for populating cubes that sum up to a given integer $k$. To illustrate it, some seed equations and populated equations are given in the *Table* **7.1**. It has been proved in the paper that every seed equation populates one more equation, that itself proves there is always more than one set of cubes that sum up to given integer. Using the populated equation as seed equation, we can further populate another equation and so on. Notwithstanding, the list given in *Table* **7.1,** there are many other integers that have been described in this paper as expressible by more than one set of three adding cubes. It is prudent to draw attention to sub paragraph **3.2a** dealing with integer of the type $3(P!/2!)^3$ and *Lemma* **3.4** that states, if integer $p = (P!/2!)^3$ and $P \geq 3$, then integer $k = 3(P!/2!)^3$ is representable by sum of cubes in as many ways as $P!$ has distinct factors of the form $d(d + 1)$. When $P \geq 5$, then $3p^3$ can be represented by adding cubes in as many sets as equal to or more than P. Also If $P!$ has n sets of distinct factors, then there will be n distinct sets of three cubes to represent integer $3(P!/2!)^3$. We have also given an example for integer $3(6!/2)^3$ which can be represented by as many as eight sets of three cubes by putting different values of $d$ using the method adopted by us. We have also proved in sub paragraph 3.7 that if $A$, $B$ and $C$ are Pythagorean's triple and $k = A^3 - B^3 + C^3$, then $k$ can also be expressed by relation $k = (3S + A)^3 + (3S - B)^3 + (-3S + C)^3$ where $S = B - (A + C)$. That also proves multiplicity of representation of $k$. Thus a given integer $k$ has multiple representation by three adding cubes.

*Table* **7.1** Some Seed Equations And Populated Equations

| $k$ | Seed Equation | Populated Equation |
|---|---|---|
| 1 | $(1)^3 + (0)^3 + (0)^3$ | $(1 - 9p^3)^3 + (3p - 9p^4)^3 + (9p^4)^3$ |
| 2 | $(1)^3 + (1)^3 + (0)^3$ | $(1 - 6p^3)^3 + (1 + 6p^3)^3 + (-6p^2)^3$ |
| 3 | $(1)^3 + (1)^3 + (1)^3$ | $(4)^3 + (4)^3 + (-5)^3$ |
| 7 | $(2)^3 + (-1)^3 + (0)^3$ | $(44)^3 + (-169)^3 + (168)^3$ |
| 8 | $(2)^3 + (0)^3 + (0)^3$ | $(-16)^3 + (-12)^3 + (18)^3$ |
| 9 | $(2)^3 + (1)^3 + (0)^3$ | $(-52)^3 + (217)^3 + (-216)^3$ |
| 10 | $(2)^3 + (1)^3 + (1)^3$ | $(4)^3 + (-3)^3 + (-3)^3$ |

## 8 Results and Conclusions

To sum up, let us have a relook what we have proved in the paper. An integer, say $X$, can be expressed by the relation, $X = a \cdot x + A$, similarly $Y$ as $b \cdot x + B$ and $Z$ as $c \cdot x + C$, where $a$, $b$, $c$, $A$, $B$ and $C$ are real rational quantities and $x$ is a variable. If $X^3 + Y^3 + Z^3 = k$, then this relation can be written in cubic equation with variable $x$. To transform this equation, first into quadratic equation, $k$ must be equal to $A^3 + B^3 + C^3$. We will, therefore, assign such values to $A, B$ and $C$ that $A^3 + B^3 + C^3 = k$. For example, integer 1 can be written as $(1)^3 + (0)^3 + (0)^3$, integer 2 can be written as $(1)^3 + (1)^3 + (0)^3$ and so on. Such equation for $k$ is called seed equation as this will populate other three cubes that sum up to $k$. To illustrate it, some seed equations are given in *Table* **1.1**. When cubic equation is transformed into quadratic, then quadratic can be transformed into linear equation by equating constant term to 0. Linear equation gives value of $x$ which can be substituted in $X = a \cdot x + A$, $Y = b \cdot x + B$ and $Z = c \cdot x + C$. Once $X, Y$ and $Z$ are known, three cubes that sum up to $k$ are known. Values of $A, B$ and $C$ depend upon the value of $k$ but values of $a, b$ and $c$ are so chosen that $X = a \cdot x + A$, $Y = b \cdot x + B$ and $Z = c \cdot x + C$ are integers.

We have also proved when $k = 3(P!/2!)^3$, then $k$ can be represented by sum of three cubes in as many ways as number of factor $d(d + 1)$. For example when $k = 3(6!/2!)^3$ then $k$ can be represented by as many as eight sets of three cubes. We also proved if $A, B$ and $C$ are Pythagorean's triple, then $k = A^3 - B^3 + C^3$, can also be represented by relation $k = (3S + A)^3 + (3S - B)^3 + (-3S + C)^3$, where $S = B - (A + C)$. Under such circumstances, an integer $k$ can be represented in multiple ways by sum of three cubes provided $k \neq 9x \pm 4$ or $k \neq 9x \pm 5$ thanks to Roger-Brown conjecture [4]. We have proved in **Sub-paragraph 3.4** why it is not possible to have three cubes that sum up to $k$ where $k = 9x \pm 4$ or $k = 9x \pm 5$.

Parametrisation for representing integer $k$ of the form $k = 9p^3 + 1$, $k = 9(p + 1)(p^2 - 1)$, $k = (9p^3 - 1)^2$, $k = 2(5 - 9p^2)$, $k = \{-10p^6 + 2\}$, $k = 2(1 + 4p^3 - 9p^4)$, $k = -\{6 + 18p^2\}$,

$k = -\{62 + 72p^2\}$, $k = 6\{11 - 12p^2\}$, $k = 1 - 2p^3$, $k = 4(2 - 5p^3)$, $k = (1 - 20p^3)$, $k = 9(q - 1)(1 + q^2) + 2$, $k = 216\{p^2 - 3\}\{2p - 3\} + 2$, $9\{2p^2 + 1\} - 3$, and $k = 2(3q + 1)^3 + 1$ as sum of three cubes, has also been given in the paper. Notwithstanding these, general parametrisation for $k = 2C^2(4C - 9p^2 \cdot q) + 2q^3$ and $k = 9\{p^2(r + q) - C\}\{p(r^2 - q^2) - C^2\} + r^3 + q^3 - C^3$ have also been determined. These can facilitate determining three cubes that sum up to a given integer $k$ of the above said forms where $p, q, r$ and $C$ are integers. We have also given a method to find three cubes when seed equation is not helpful or is difficult to find out. In addition to these, we have given exhaustive examples to illustrate the novel method of finding three cubes that sum up to a given number $k$.

## 9  Acknowledgements

We acknowledge the help provided by the scientific calculators at websites https://www.desmos.com and https://keisan.casio.com/calculator in calculating the values of tedious and large exponential and factorial terms. We further acknowledge the guidance and help provided by the learned Reviewer and learned Editor in bringing the paper to its present form.

## 10  References

[1]   A. Avagyan, G. Dallakyan, A New Method in the Problem of Three Cubes, *Universal Journal of Computational Mathematics,* **5 (3)**, (2017), 45-56.

[2]   A. S. Verebrusov, Объ уравненіи $x^3 + y^3 + z^3 = 2u^3$, [On the equation $x^3 + y^3 + z^3 = 2u^3$], *Matematicheskii Shormik*, (in Russian), **26 (4)**, (1908), 622–624.

[3]   A. R. Booker, A. V. Sutherland, On a question of Mordell, *arXiv:2007.01209*, (2020).

[4]   D. R. Heath-Brown, W. M. Lioen, and H. J. J. Te Riele, On Solving the Diophantine Equation $x^3 + y^3 + z^3 = k$ on a Vector Computer, *Mathematics of Computation,* **61,** (1993), 235-244.

[5]   A. Elsenhans, J. Jahnel, New sums of three cubes, *Mathematics of Computation,* **78(266),** (2009), 1227–1230.

[6]   A. Georgiou, The untracked problem with 33: Mathematicians solves 64-year-old Diophantine puzzle, *Newsweek*, (April 3, 2019),

[7]   R. Houston, 42 is the answer to the question 'what is $(-80538738812075974)^3 + (80435758145817515)^3 + (12602123297335631)^3$ ? *The Aperiodical,* (September 6, 2019).

[8]   L. J. Mordell, On sums of three cubes", *Journal of the London Mathematical Society*, **17 (3)**, (1942), 139–144.

[9]   D. Lu, Mathematics cracks centuries-old problem about the number 33, *New Scientist,* (March 14, 2019).

[10]  K. Mahler, Note on Hypothesis K of Hardy and Littlewood, *Journal of the London Mathematical Mathematical Society*, **11(2),** (1936), 136–138.

[11]  J. Pavlus, Sum of Three Cubes Problem Solved for 'Stubborn Number 33', *Quanta Magazine,* (March 10, 2019).